\RequirePackage{ifpdf}
\ifpdf 
\documentclass[pdftex]{sigma}
\else
\documentclass{sigma}
\fi

\begin{document}

\allowdisplaybreaks

\renewcommand{\PaperNumber}{081}

\FirstPageHeading

\renewcommand{\thefootnote}{$\star$}

\ShortArticleName{Future Directions of Research in Geometry}

\ArticleName{Future Directions of Research in Geometry:\\
A Summary of the Panel
Discussion\\ at the
2007 Midwest Geometry Conference\footnote{This paper is a
contribution to the Proceedings of the 2007 Midwest
Geometry Conference in honor of Thomas~P.\ Branson. The full collection is available at
\href{http://www.emis.de/journals/SIGMA/MGC2007.html}{http://www.emis.de/journals/SIGMA/MGC2007.html}}}

\Author{Edited by Lawrence J. PETERSON}

\AuthorNameForHeading{L.J. Peterson}

\Address{University of North Dakota, Grand Forks, North
  Dakota, USA} 
\Email{\href{mailto:lawrence_peterson@und.nodak.edu}{lawrence\_peterson@und.nodak.edu}}

\URLaddress{\url{http://www.und.nodak.edu/instruct/lapeters}}

\ArticleDates{Received August 09, 2007; Published online August 15, 2007}

\Abstract{The 2007 Midwest Geometry Conference included a panel discussion
devoted to open problems and the general direction of future research
in f\/ields related to the main themes of the conference.  This paper
summarizes the comments made during the panel discussion.}


\Classification{53A35}

\begin{flushright}
\it Dedicated to the memory of Thomas P.~Branson
\end{flushright}

The 2007 Midwest Geometry Conference was held at the University of
Iowa in Iowa City, Iowa, USA (May 18--20, 2007).  In addition to the
usual talks given at such meetings, the conference also included a
panel discussion devoted to open problems and the general direction of
future research in geometry and related f\/ields.  This paper
is an edited summary of the panel discussion and the contributions of
the members of the audience.

The panel included Sun-Yung Alice Chang (Princeton), Michael Eastwood
(Adelaide), A.~Rod Gover (Auckland), C. Robin Graham (University of
Washington), Rafe Mazzeo (Stanford), Bent {\O}rsted (Aarhus), and Paul
Yang (Princeton).  The audience members who spoke were James G. Branson
(University of California, San Diego), Stephen Fulling (Texas A \& M),
Colin Guillarmou (Nice), Claude LeBrun (Stony Brook), Andrew Waldron
(University of California, Davis), and Shihshu Walter Wei (Oklahoma).
Rafe Mazzeo chaired the panel.

Most of the panel discussion related to the main theme of the 2007
Midwest Geometry Conference, namely the work of Thomas P. Branson
(1953--2006).  The panel and the audience discussed individual open
problems and topics that people will likely investigate in the years
to come.  An edited summary of this discussion follows.  We provide
the various comments in approximately the same order that they
occurred during the actual group discussion.  The discussion began
with some remarks by Rafe Mazzeo.

\textbf{Rafe Mazzeo:} Three prominent topics have come up in the
course of the 2007 Midwest Geometry Conference so far, namely
\begin{itemize}\vspace{-2mm}\itemsep=0pt
\item Branson's $Q$-curvature;
\item Poincar\'{e}--Einstein metrics;
\item Bach-f\/lat metrics.
\end{itemize}

\textbf{James Branson:}  What does all of this have to do with
physics?

\textbf{\textit{Editor's note:}} For a partial answer to this
question, see the replies by Andrew Waldron and Stanley Deser near the
end of this paper.  After the question by Jim Branson, the
discussion turned to some specif\/ic open problems.

\textbf{Claude LeBrun:}  Are there compact four-dimensional manifolds
which are Bach-f\/lat and neither
\begin{enumerate}\itemsep=0pt\vspace{-2mm}
\item[1)] locally conformally Einstein nor
\item[2)] half f\/lat, i.e.\ self-dual or anti-self-dual?\vspace{-2mm}
\end{enumerate}
It might be reasonable to start by investigating the toric case,
that is, by considering only those metrics which admit an isometric
$U(1)\times U(1)$ action.

Another open problem is to classify asymptotically locally Euclidean
scalar-f\/lat K\"{a}hler metrics on complex surfaces.  Spaces of this type
turn out to naturally arise as blow-up limits in many problems related
to the convergence of sequences of Bach-f\/lat metrics, and a complete
classif\/ication would provide a powerful tool for proving that
curvature can only accumulate and bubble of\/f in very specif\/ic
circumstances.

\textbf{Alice Chang:} A very general question is to ask ``What is the
geometric content of $Q$-curvature?''  For example, we know that one can
associate the scalar curvature with the conformally invariant constant
called the ``Yamabe constant''.  When this constant is positive, it
describes the best constant (in a conformally invariant sense) of the
Sobolev embedding of $W^{1,2}$ into $L^{2n/(n-2)}$ space; this in
itself can be viewed as a $W^{1,2}$ version of the isoperimetric
inequali\-ty.  It would be interesting to know if $Q$-curvature, or the
conformally invariant quantity $\int Q$ associated with it, satisf\/ies
some similar inequa\-li\-ties with geometric content.

\textbf{Mike Eastwood:} It would be good to have a characterization of
$Q$-curvature.  Although we know that $Q$-curvature has many properties,
we do not understand enough of its properties to characterize it.

\textbf{Alice Chang:} It would be nice to have a uniqueness result
which would characterize $Q$-curvature.

\textbf{Mike Eastwood:} One can use the Fef\/ferman--Graham ambient
metric to canonically construct Branson's $Q$-curvature.  Why is this
so?

\textbf{Alice Chang:} Can one describe the ``renormalized volume'' in
terms of the ambient metric?

\textbf{Robin Graham:} Concerning Mike Eastwood's request for a
characterization of $Q$-curvature, I made a conjecture of such a
characterization at the workshop on ``Conformal structure in geometry,
analysis, and physics'' at the American Institute of Mathematics in
2003.  See problems 10 and 11 of the open problems list available at
\begin{quotation}
  \url{http://www.aimath.org/pastworkshops/confstruct.html}
\end{quotation}
To my knowledge, this conjectured characterization is still open and,
if true, would provide a~good characterization of $Q$-curvature.  The
conjecture comes in two parts.  The f\/irst part is a~conjectured
characterization of the critical ``GJMS'' operator, which appears in the
conformal transformation law for $Q$-curvature.  The GJMS construction
shows that the coef\/f\/icients of each GJMS operator can be written just
in terms of the Ricci tensor and its covariant derivatives, and that
the $Q$-curvature can likewise be written just in terms of the Ricci
tensor and its covariant derivatives.  My specif\/ic two-part conjecture
is as follows.
\begin{enumerate}\itemsep=0pt\vspace{-2mm}
\item Let $n\geq 4$ be even.  Suppose that $P$ is a scalar natural
    operator in dimension $n$ whose principal part is $\Delta^{n/2}$.
    Suppose also that $P$ satisf\/ies the same conformal covariance
    relation as the critical GJMS operator.  If all of the
    coef\/f\/icients of $P$ can be written in terms of the Ricci tensor
    and its covariant derivatives, then $P$ is the critical GJMS
    operator.
\item A local scalar conformal invariant of weight $-n$ which can
    be written only in terms of the Ricci tensor and its covariant
    derivatives must vanish.\vspace{-2mm}
\end{enumerate}
If true, the f\/irst conjecture would uniquely characterize the critical
GJMS operator, and then the second would imply that $Q$-curvature is
uniquely characterized by its conformal transformation law involving
this critical GJMS operator.  Both conjectures are true in dimensions
4 and 6.

Notes: The GJMS operators are the conformally covariant partial
dif\/ferential operators of Graham, Jenne, Mason, and Sparling
\cite{GJMS}.  The critical GJMS operator is the operator of order~$n$
in even dimension $n$.

\textbf{Claude LeBrun:} Another open problem is to characterize those
K\"{a}hler metrics for which the Fef\/ferman--Graham obstruction tensor
vanishes.  The goal would be to do this for arbitrary dimensions.  In
real dimension four, the Fef\/ferman--Graham obstruction tensor is just
the Bach tensor, and in \cite{Derd}, Andrzej Derdzinski showed that it
vanishes for a K\"{a}hler metric if and only if the metric is
conformally Einstein on the open set where the scalar curvature is
nonzero. This observation plays a central role in my joint work with
Chen and Weber \cite{LeBrun}, as I described in my lecture yesterday.
Since the Fef\/ferman--Graham obstruction tensor exactly represents the
gradient of the total $Q$-curvature, and since many dif\/ferential
geometric problems simplify when restricted to the K\"{a}hler arena,
this might be an excellent context in which to obtain a~more concrete
understanding of these higher-order curvature objects, which we now
but glimpse through a glass, darkly.

\textbf{Rod Gover:} An emerging theme in my more recent work with Tom
Branson was the construction of a new class of dif\/ferential complexes
that we termed ``detour complexes''.  These complexes exist on, for
example, even-dimensional conformally f\/lat manifolds and are
conformally invariant.  They are related to (but distinct from)
Bernstein--Gelfand--Gelfand (BGG) complexes.  In Riemannian signature,
the cases treated yield elliptic complexes.  One open problem is to
show that all conformal detour complexes are elliptic.  So far there
are constructions of several families of curved analogues.  The
general theory of curved analogues should be developed.

Related to the de Rham detour complexes is a generalization of
$Q$-curvature.  In my joint work with Tom, we obtained what are termed
``$Q$-operators''.  These exist on even-dimensional conformal manifolds;
they act on closed forms and determine conformally invariant
cohomology pairings.  Any progress in understanding $Q$-operators would
be extremely interesting.  In the case of zero-forms, the $Q$-operators
recover the usual $Q$-curvature, and the cohomology pairing in that case
is just the integral of the $Q$-curvature.  On conformally f\/lat spaces
this recovers the Euler characteristic.  So an important question is
whether on conformally f\/lat structures the $Q$-operators have a link to
topology which generalizes this.

The $Q$-operators yield problems which generalize the problem of
prescribing constant $Q$-curvature.  In these problems, one begins with
a conformal class of metrics.  One then seeks a~metric in this
conformal class and a closed dif\/ferential $k$-form which, for this
particular metric, is an eigenform for the $Q$-operator on $k$-forms.  The
eigenvalue is the analogue of constant $Q$.  It would be interesting to
know if this yields a manageable PDE problem.

\textbf{Bent {\O}rsted:}
In the years to come there will be further research in Sobolev
inequalities.  People will study strong Hardy--Littlewood--Sobolev
inequalities and their analogues in the context of Cauchy--Riemann (CR)
geometry and other parabolic geometries.  One goal will be to maximize
or minimize linear functionals.

Tom Branson was interested in calculating determinants of matrices in
the setting of intertwining operators; when the $K$-types have higher
multiplicity, an intertwining operator is not just a number for each
$K$-type, but rather a matrix, and one would like to compute its
determinant.  This has been studied by Wallach, Cohn, and Vogan.  See
\cite{VoWall}, for example.  There is some hope that the method of spectrum
generating could give results in this direction.  Actually, Branson
had some computations for a multiplicity-two example involving
dif\/ferential forms.

\textbf{Paul Yang:} It would be interesting to study $Q$-curvature in
odd dimensions.  The plan would be to def\/ine the GJMS operator as a
nonlocal operator in the context of Poincar\'{e}--Einstein metrics as
def\/ined by Fef\/ferman and Graham.  This should provide interesting
analytic questions.

\textbf{Robin Graham:} Spyros Alexakis reports that he has
characterized scalar Riemannian invariants whose integral over a
compact manifold is a conformal invariant.  He has characterized them
as a linear combination of a divergence, a multiple of the Pfaf\/f\/ian,
and a local conformal invariant of weight $-n$, where $n$ is the
dimension.  It is an interesting question whether such a decomposition
is unique.  It is always unique in dimensions 4 and 6.

\textbf{Rod Gover:} For suitable conformally invariant dif\/ferential
operators and for operators or complexes with appropriate ``nearly
invariant'' properties, there are determinants (or torsions) which are
not conformally invariant but nevertheless have useful ``Polyakov
type'' formulas for their conformal variation.  In my joint work with
Tom Branson, Tom has proposed using invariant theory and related tools
to investigate the existence of a universal form for such formulas.

\textbf{Robin Graham:} Neil Seshadri has def\/ined a torsion in the
contact (CR) case analogous to the conformal half-torsion.  See
\cite{Sesh}.  Seshadri's torsion also has a local variation under
change of contact form and almost complex structure, which is given in
terms of heat kernel coef\/f\/icients.  Further work needs to be done to
make this more explicit and to better understand such formulas.

\textbf{Bent {\O}rsted:} In addition to the above problems, people
should study other invariants as well.

\textbf{Stephen Fulling:} Ten to twenty years ago, many people,
including Tom Branson, devoted a great deal of attention to
calculating heat kernel coef\/f\/icients explicitly to the highest order
feasible.  About ten years ago we all seemed to get tired of it, but
it may now be time to look at the subject again.  It seemed to me that
we were right on the verge of reaching an understanding of the
structure at some more profound algebraic or combinatorial level.

In the case of a Schr\"{o}dinger operator with just a potential (in
f\/lat space) one can write down an explicit combinatorial formula for
the $n$'th coef\/f\/icient, even before passing to the diagonal and taking
the trace.  For a Riemannian manifold, a computer can in principle
grind out the formulas recursively, but the result is not
illuminating.  The problem arises in reducing the formula to a
canonical form.  There is an accumulation of terms of high degree in
the Riemann tensor coming from the commutation of covariant
derivatives in the process of making the low-degree terms linearly
independent; at the end, the coef\/f\/icient of such a term (a sum of many
disparate contributions) conveys no clear information.  Yet there
could still be something to learn here.

\textbf{Shihshu Walter Wei:} My question is related to Professor Alice
Chang's question: Is there any connection between $Q$-curvature and
$p$-harmonic geometry?  For example, let us recall that a manifold $M$
is $p$-hyperbolic, $1 \le p < \infty$, if the $p$-capacity of any
compact subset is positive (or equivalently, if $M$ supports a
nonconstant positive $p$-superharmonic function) and $p$-parabolic
otherwise.  If $p=n=\dim M$, then $n$-parabolicity of $M$ is a
(quasi)conformal invariant (cf. \cite{T}).  We know that if the
sectional curvature of a complete simply-connected manifold $M$ is
bounded above by a negative constant, then $M$ supports an
isoperimetric inequality, and via isoperimetric prof\/ile we see that
the manifold is $p$-hyperbolic for every $p \ge 1$
(cf.\ e.g.~\cite{CW}). On the other hand, we also know that the
Brunn--Minkowski inequality in the theory of convex bodies is connected
with other fundamental inequalities, such as the isoperimetric,
Sobolev, and Pr\'{e}kopa--Leindler inequalities, and can be extended to
the $p$-capacity of convex bodies, $1 < p < n$ (cf.~\cite{CS}).

It would be interesting to know if $Q$-curvature, or the conformally
invariant quantity $\int_M Q$ associated with it, is linked to
$p$-parabolicity, $p$-hyperbolicity, $p$-capacity, $p$-energy, or
other functionals arising in the context of calculus of variations,
the theory of PDEs, or geometric measure theory, or satisf\/ies some
inequality of geometric signif\/icance.

\textbf{Colin Guillarmou:} Another interesting problem is to study the
dimension of the kernel of the critical GJMS or Branson--Gover
operators.  The dimension of the kernel of the GJMS operators is a
conformal invariant.  It appears in the Selberg trace formula for
(non-compact) convex co-compact hyperbolic manifolds. These manifolds
have (locally conformally f\/lat) conformal inf\/inity.  Thus the natural
GJMS operators and the dimension of their kernels show up in the trace
formula (see \cite{G} and \cite{PP}).  Note also that the kernel of
the critical GJMS operator is a conformally invariant space.  The same
problem for the Branson--Gover operators on forms is also a parallel
thing to consider.

It is also of interest to study the determinant of the GJMS operators.
Branson, Chang and Yang worked on the Yamabe and square of Dirac
cases.  It could be worthwhile to look at higher-order GJMS operators
and see if it is possible to obtain results in the same spirit as in
the case of the determinant of the Laplacian. For instance, it is
possible to compute the determinant of the GJMS operators explicitly
for a class of odd-dimensional locally conformally f\/lat manifolds.
(See \cite{G}.)

\textbf{Robin Graham:} Over the last few years, there have been a
number of constructions of global invariants of contact manifolds
which arise (or are expected to arise) as integrals of local
invariants of a choice of contact form and compatible almost complex
structure.  For example, Hirachi, Boutet de Monvel, Ponge, Seshadri,
and Biquard--Herzlich--Rumin have developed constructions of such global
invariants.  None of these global invariants is known to be nonzero,
and some of them have been shown to vanish.

This gives rise to a natural conjecture about such invariants, which
is formulated explicitly by Seshadri in Appendix A of \cite{Sesh}.  The
conjecture is motivated by a result of Gilkey which goes as follows.
Consider any local scalar Riemannian invariant whose integral over a
compact manifold is independent of the metric.  This invariant is
necessarily a divergence plus a multiple of the Pfaf\/f\/ian; the integral
must therefore be a multiple of the Euler characteristic.

The analogous conjecture in the contact case reads as follows.
Suppose that one has a local scalar invariant of a contact manifold
together with a choice of contact form and compatible almost complex
structure whose integral over a compact manifold is a contact
invariant.  Then the local invariant is a divergence, and hence the
integral must vanish.

A positive answer to this conjecture would indicate a unif\/ied approach
to showing that all the invariants which arise in the above-mentioned
constructions necessarily vanish.  It would also show that there is no
contact analogue of the Pfaf\/f\/ian.

\textbf{Claude LeBrun:} It seems worth pointing out that the results
Robin just alluded to only hold if one requires the global invariant
to be independent of the orientation of the manifold; in particular,
the hypotheses exclude things like Pontrjagin numbers.  Perhaps one
should try to classify invariants for which the invariant is allowed
to change sign when the orientation of the manifold is reversed.
Next, one might ask what happens when the dependence on orientation is
allowed to be arbitrary.

\textbf{Andrew Waldron:} It would be interesting to study Lorentzian
versions of many of these questions.  Our universe has Lorentzian
signature.  Many formulas look exactly the same in Lorentzian and
Riemannian signatures; other aspects dif\/fer signif\/icantly, however,
especially when one wants to solve, rather than just write down,
interesting systems of PDEs on curved manifolds.

\textbf{Rod Gover:} It would be interesting to consider the above
questions in the context of special structures such as almost Einstein
metrics.

\textbf{\textit{Editor's note:}} As we noted above, Jim Branson posed
the question ``What does all of this have to do with physics?''  In
the weeks following the 2007 Midwest Geometry Conference, Stanley
Deser, of Brandeis University, has provided us with a short essay
which comments on ways in which Tom Branson's work relates to physics.
Andrew Waldron has also provided additional comments.  Waldron's
comments and Deser's essay follow.

\textbf{Andrew Waldron:} There is a deep relationship between
conformal geometry on the one hand and the AdS/CFT correspondence in
physics on the other.  The latter, f\/irst suggested by Juan Maldacena,
relates gravity/string theories in an $(n+1)$-dimensional
asymptotically Anti de Sitter background to conformal f\/ield theories
on their boundaries. In particular the $Q$-curvature and Poincar\'{e}
metric are related to conformal anomalies of the boundary conformal
f\/ield theories.

\textbf{\textit{Note:}}  AdS/CFT means ``Anti de Sitter/conformal f\/ield theory''.

\medskip

\textbf{Stanley Deser's Essay:} Tom Branson~--~An Appreciation from Physics. 

Theoretical physicists are supposed to know, at least vaguely, some of
the mathematical riches available to help our research. The reality is
often rather dif\/ferent. First you try to locate someone in that
distant country who can help you, if only you can f\/igure out how to
frame your vague queries intelligibly in his or her language.  And
then, if you are lucky, you just have to translate back the usually
far too technical answer.

Tom Branson was the exception I was lucky enough to encounter early
enough in a project on conformal anomalies, a subject of great
physical interest. True, he did at f\/irst refer me to his papers or
quote theorems~-- all interaction was by email or (rarely) by phone~-- but
he wisely realized that I needed words of one syllable explanations,
and he cheerfully provided them.

I will be monosyllabic myself in describing the context; a paper that
provides technical and citation details is \cite{Deser}.  Conformal
anomalies are described by ef\/fective actions in (even) dimension~$D$
that hinge on very special operators such as the inverse Laplacian in
$D=2$ and the Paneitz operator~-- the self-adjoint conformal completion
of the inverse square Laplacian~-- in $D=4$.  Both Tom and Paneitz (who
also died tragically, already as a graduate student) were students of
Irving Segal at MIT, where Irving's passion for all things conformal
inspired them. Actually, there are two dif\/ferent anomaly types with
very dif\/ferent properties and origins. The Paneitz operator and its
generalizations (to which Tom also contributed) arise in type A.  Type
B requires a very dif\/ferent approach, and only starts at $D=4$. Here
the desired operators are not only bilocal but heavily tensorial,
unlike the scalar ones of type A.

This is where Tom really saved the day.  It was morally certain to a
physicist that there had to exist a fourth derivative order (in $D=4$)
of tensor rank eight which when acting to any power on a tensor with
the symmetries of the conformal curvature (Weyl) left it conformally
invariant.  It was Tom who actually had the desired object (and of
course its generalizations for arbitrary $D$) and who told me that he
and Rod Gover were going to publish it.  Talk of mining mathematical
treasure with no ef\/fort! This work enabled me to to give the closed
form solution to the type B ef\/fective action, whose original discovery
in 1976 was only given as a hand-waving scalar approximation to the
real thing.

Emboldened by these riches, I was pumping Tom for many other important
unf\/illed physics needs that would have propelled us, when the terrible
news of his death came. Despite the action-at-a-distance nature of our
relation, the loss was felt as one of an old friend!

\subsection*{Acknowledgements}

This paper is based in part on notes the editor took during the panel
discussion at the 2007 Midwest Geometry Conference.  The editor would
like to thank Alice Chang, Stanley Deser, Mike Eastwood, Stephen
Fulling, Rod Gover, Robin Graham, Colin Guillarmou, Palle Jorgensen,
Claude LeBrun, Rafe Mazzeo, Bent {\O}rsted, Andrew Waldron, Shihshu
Walter Wei, and Paul Yang for the additional notes and clarif\/ication
that they provided in the weeks following the conference and panel
discussion.

This material is based upon work supported by the National Science
Foundation under grants No.\ DMS-0509068 and No.\ PHY04-01667.

Any opinions, f\/indings, and conclusions or recommendations expressed
in this material are those of the authors and do not necessarily
ref\/lect the views of the National Science Foundation.

\pdfbookmark[1]{References}{ref}
\LastPageEnding

\end{document}